\documentclass[11pt]{article}
\usepackage{amssymb, amsfonts, amsmath }
\usepackage{latexsym}

\providecommand{\MR}{\relax\ifhmode\unskip\space\fi MR }

\newtheorem{thm}{Theorem}[section]
\newtheorem{prop}[thm]{Proposition}

\newtheorem{cor}[thm]{Corollary}
\newtheorem{example}[thm]{Example}

\newtheorem{dfn}[thm]{Definition}

\newcommand{\demo}{ \noindent {\it   Proof. }}
\newcommand{\qed}{$\Box$\bigskip
}
\newcommand{\au}{\mbox{Aut}\,}
\newcommand{\en}{\mbox{End}\,}
\newcommand{\auH}{\mbox{Aut}_H}
\newcommand{\enH}{\mbox{End}_H}

\newcommand{\fix}{\mbox{Fix\,}}
\newcommand{\acl}{\mbox{\it a-Cl\,}}
\newcommand{\ecl}{\mbox{\it e-Cl\,}}
\newcommand{\leqff}{\leqslant_{\rm ff}}
\newcommand{\leqalg}{\leqslant_{\rm alg}}

\title{Computing fixed closures in free groups}

\author{E.\ Ventura \thanks{The author gratefully acknowledges partial support from the MEC (Spain) and the EFRD (EC)
through project number MTM2008-01550.} \\ \small{Dept.\ Matem\`{a}tica Aplicada III} \\ \small{Universitat Polit\`{e}cnica de
Catalunya,} \\ \small{Barcelona, Catalonia.} \\ \small{e-mail: enric.ventura@upc.edu} }

\date{\today}

\begin{document}

\maketitle


\begin{abstract}
Let $F$ be a finitely generated free group. We present an algorithm such that, given a subgroup $H\leqslant F$, decides
whether $H$ is the fixed subgroup of some family of automorphisms, or family of endomorphisms of $F$ and, in the
affirmative case, finds such a family. The algorithm combines both combinatorial and geometric methods.
\end{abstract}

\section{Introduction}\label{intro}

For all the paper, let $A=\{ a_1,\ldots ,a_n \}$ be an alphabet with $n$ different letters, and $F$ be the free group
(of rank $r(F)=n$) with basis $A$.

Let $\en (F)$ denote the endomorphism monoid of $F$, and $\au (F)$ the automorphism group of $F$, so $\au (F)$ is the
group of units of $\en (F)$. Throughout, we let elements of $\en (F)$ act on the right on $F$, so $x \mapsto (x)\phi$.
Accordingly, compositions are like $(x)\phi\psi =(x\phi)\psi$.

In the last decade a lot of literature has appeared studying the fixed subgroup of a single, or a family, of
automorphisms, or endomorphisms, of $F$ (see the survey~\cite{V} for details). But very few algorithmic results are
known in this direction. Only few years ago, O. Maslakova (see~\cite{Mas}) has found an algorithm to compute a set of
generators for the fixed subgroup of an automorphism of $F$ (which is quite complicated, and whose complexity is quite
high). One can easily extend this to compute generators for the fixed subgroup of a finite family of automorphisms, but
the related questions on endomorphisms are still open.

In this note, we shall address the dual problem. We present an algorithm such that, given a subgroup $H\leqslant F$,
decides whether $H$ is the fixed subgroup of some finite family of automorphisms, or finite family of endomorphisms of
$F$ and, in the affirmative case, finds such a family. We advice the reader that the provided algorithms are
theoretical and far from effective, in the sense that their complexities will be way too high to think about possible
effective implementations. A natural open question is whether there exist more natural and efficient algorithms, say
polynomial time, for solving such problems.

Recognizing whether $H$ is the fixed subgroup of a family of automorphisms is not difficult, because a classical result
by McCool already established that the stabilizer of a finitely generated subgroup of $F$ is finitely generated and
computable (as subgroup of $\au (F)$), see Theorem~\ref{McCool} below. However, recognizing whether $H$ is the fixed
subgroup of a family of endomorphisms is much trickier, because the stabilizer of $H$, in general, need not be finitely
generated as a submonoid of $\en (F)$ (see Example~\ref{ex} below).

More precisely, the two algorithms given in this paper take a finitely generated subgroup $H\leqslant F$ as the input,
and compute a basis of its automorphism (resp. endomorphism) closure, i.e. the smallest subgroup $K$ such that
$H\leqslant K\leqslant F$ and $K=\fix (S)$ for some $S\subseteq \au F$ (resp. $S\subseteq \en F$), see the precise
definitions below. The main technique used to deal with these problems is the graphical tool called ``fringe of a
subgroup", which allows to compute the collection of algebraic extensions of a given subgroup $H$.

In Section~\ref{tools} we define the concepts, and state the results that will be used, with the corresponding
references. Along Section~\ref{algo} we prove the main result and give the two announced algorithms. Finally, in
Section~\ref{p} we collect a list of related questions and open problems.

\section{Needed tools}\label{tools}

\begin{dfn} {\rm
For any $S \subseteq \en (F)$, let $\fix (S)$ denote the subset consisting of those elements of $F$ which are fixed by
every element of $S$ (read $\fix (S)=F$ for the case where $S$ is empty). Then $\fix (S)$ is a subgroup of $F$, called
the \emph{fixed subgroup of $S$}.

A subgroup $H$ of $F$ is called an {\it endo-fixed} subgroup of $F$ if $H =\fix (S)$ for some $S\subseteq \en (F)$. If
$S$ can be chosen to lie in $\au (F)$ we further say that $H$ is an {\it auto-fixed} subgroup of $F$.

A subgroup $H$ of $F$ is called a {\it $1$-endo-fixed} subgroup of $F$ if $H = \mbox{Fix}(\phi)$ for some $\phi\in \en
(F)$ (here, and throughout, to simplify notation we write $\mbox{Fix}(\phi)$ rather than $\mbox{Fix}(\{\phi\})$). If
$\phi$ can be chosen to lie in $\au (F)$, we further say that $H$ is a {\it $1$-auto-fixed} subgroup of $F$. For
example, any maximal cyclic subgroup of $F$ is $1$-auto-fixed, since it is the subgroup fixed by a suitable inner
automorphism. And non-maximal cyclic subgroups of $F$ are not even endo-fixed, because every endomorphism fixing a
power of an element must fix the element itself.

Notice that, since $\fix(S)=\cap_{\alpha \in S} \fix(\alpha)$, an auto-fixed (resp. endo-fixed) subgroup is an
intersection of $1$-auto-fixed (resp. 1-endo-fixed) subgroups, and vice-versa. And, clearly, the families of auto-fixed
and endo-fixed subgroups of $F$ are closed under arbitrary intersections. $\Box$}
\end{dfn}

A natural question that arises in this context asks about the relation between the four mentioned families of subgroups
of $F$, namely $1$-auto-fixed, $1$-endo-fixed, auto-fixed and endo-fixed subgroups. Apart from the obvious inclusions,
the relationship among these families is partially known, though not completely since there still are interesting
questions in this direction that remain open. For example, it is not known (and conjectured) whether the families of
$1$-auto-fixed and auto-fixed subgroups (resp. $1$-endo-fixed and endo-fixed subgroups) do coincide; in other words, it
is not known whether the family of $1$-auto-fixed (resp. $1$-endo-fixed) subgroups is closed under intersections. As
far as we are aware, this is only known to be true when the ambient rank is $n=2$, and when the involved fixed
subgroups have maximal rank, see~\cite{MV2}. In this direction, A.~Martino and E.~Ventura showed in~\cite{MV2} that
every auto-fixed (resp. endo-fixed) subgroup of $F$ is a free factor of a $1$-auto-fixed (resp. $1$-endo-fixed)
subgroup of $F$. However, they also gave an example of a free factor of a $1$-auto-fixed subgroup which is not even
endo-fixed.

We point out that, in the definitions of auto-fixed and endo-fixed subgroups, one can always assume that the involved
families of morphisms are finite. This was proven by A.~Martino and E.~Ventura in~\cite[Corollary~4.2]{MV}, answering a
question previously posed by G. Levitt. So, from now on, a ``family" of endomorphisms will always mean a ``finite
family":

\begin{prop}[Martino-Ventura,~\cite{MV}]\label{2n}
Let $F$ be a finitely generated free group. For every $S\subseteq \en (F)$ there exists a finite subset $S_0\subseteq
S$ with $|S_0|\leqslant 2r(F)$, and such that $\fix (S_0) =\fix (S)$.
\end{prop}

We do not include in this discussion the families of $1$-mono-fixed and mono-fixed subgroups, because they are known to
coincide with the families of $1$-auto-fixed and auto-fixed subgroups, respectively (see~\cite[Theorem~11]{MV3}).

On the other hand, it is known that the families of $1$-endo-fixed and $1$-auto-fixed subgroups (and the families of
endo-fixed and auto-fixed subgroups, as well) \emph{do not} coincide: in~\cite{MV3}, the authors exhibited the first
known examples of $1$-endo-fixed subgroups which are not $1$-auto-fixed; see also~\cite{CD} for more interesting
calculations about this phenomena. Hence, determining whether a given subgroup $H$ is an auto-fixed or an endo-fixed
subgroup are two different algorithmic problems (the first being much simpler than the second, as will be seen below).

The deepest and most important result about $1$-auto-fixed subgroups in the literature was obtained by M.\ Bestvina and
M.\ Handel in~\cite{BH}, where they developed the theory of train tracks for graphs, and showed that every
$1$-auto-fixed subgroup of $F$ has rank at most $r(F)$, which had previously been conjectured by G.\ P.\ Scott. Soon
after the announcement of this result, and using it, W. Imrich and E. Turner showed, in~\cite{IT}, that any
$1$-endo-fixed subgroup of $F$ also has rank at most $r(F)$. Later, W. Dicks and E. Ventura in~\cite{DV}, using the
techniques of~\cite{BH}, showed that any auto-fixed subgroup of $F$ has rank at most $r(F)$; in fact, they proved a
stronger result, namely that any mono-fixed subgroup of $F$ is $F$-inert (a subgroup $H\leqslant F$ is \emph{$F$-inert}
if $r(H\cap K)\leqslant r(K)$ for every $K\leqslant F$). And after this, G.\ M.\ Bergman~\cite{B}, using the result
of~\cite{DV}, showed that any endo-fixed subgroup of $F$ also has rank at most $r(F)$ (however, it is not known whether
endo-fixed subgroups of $F$ are necessarily $F$-inert; it is conjectured to be so in the \emph{inertia conjecture},
see~\cite{MV} and~\cite{V}). This brief history is appropriate for our purposes, but is far from complete; for example,
it does not mention the ground-breaking work of S.\ M.\ Gersten, who first showed that $1$-auto-fixed subgroups are
finitely generated.

As we mentioned in the introduction, few algorithmic results are known about fixed subgroups of free groups. The main
one is the computability of fixed subgroups of automorphisms which, by the moment, it has only theoretical interest
because no precise bound on the complexity is known, and one expects it to be quite high. The corresponding fact for
endomorphisms is still an open problem.

\begin{thm}[Maslakova, \cite{Mas}]\label{mas}
Let $\varphi \colon F \to F$ be an automorphism of a finitely generated free group $F$. Then, a basis for $\fix
(\varphi)$ is computable.
\end{thm}

An interesting notion to study these questions is the notion of ``closure" of a subgroup.

\begin{dfn} {\rm
Let $H\leqslant F$. We denote by $\auH (F)$ the subgroup of $\au (F)$ consisting of all automorphisms of $F$ which fix
$H$ pointwise,
 $$
\auH (F)=\{ \varphi \in \au (F) \mid H\leqslant \fix (\varphi) \},
 $$
usually called the \emph{stabilizer} of $H$. Analogously, we denote by $\enH (F)$ the submonoid of $\en (F)$ consisting
of all endomorphisms of $F$ which fix every element of $H$. Clearly, $\auH (F) \leqslant \enH (F)$.

Now, $\au_{(-)} (F )$ is a function from the set of subgroups of $F$ to the set of subsets of $\au (F)$, and $\fix(-)$
is a function in the reverse direction. This pair of functions form a Galois connection, and their images are called
\emph{closed} subsets (in $\au (F)$ and $F$, respectively). Clearly, $\au (F)$-closed subgroups of $F$ are precisely
the auto-fixed subgroups. Mimicking the classical Galois notions, we define the \emph{auto-closure} of $H$ in $F$,
denoted $\acl_{F}(H)$, as $\fix(\auH (F))$, i.e. the smallest auto-fixed subgroup of $F$ containing $H$.

Replacing $\au$ to $\en$ everywhere in the previous paragraph we obtain another Galois connection, and we similarly
define the \emph{endo-closure} of $H$ in $F$, denoted $\ecl_{F}(H)$, as $\fix(\enH (F))$, i.e. the smallest endo-fixed
subgroup of $F$ containing $H$. Since $\auH (F)\leqslant \enH (F)$, an obvious relation between closures is that
 $$
\ecl_{F}(H) =\fix(\enH (F)) \leqslant \fix(\auH (F)) =\acl_{F} (H).
 $$
However, the equality does not hold in general, because of the existence of $1$-endo-fixed subgroups which are not
auto-fixed.

Note that, by the results mentioned above, the ranks of $\acl_{F}(H)$ and $\ecl_{F}(H)$ are always less than or equal
$r(F)$, even if that of $H$ is not. $\Box$ }
\end{dfn}

The main goal of this note is to show that, for any finitely generated $H\leqslant F$ (given by a set of generators),
one can algorithmically compute a basis for both $\acl_{F}(H)$ and $\ecl_{F}(H)$. Using this algorithm, one can
immediately decide whether the given $H$ is auto-fixed (resp. endo-fixed) or not: $H$ is auto-fixed (resp. endo-fixed)
if and only if $\acl_{F}(H)=H$ (resp. $\ecl_{F}(H)=H$).

To do this, we need to use the concepts of retract and stable image, and the graphical technique called ``fringe of a
subgroup" to compute the set of algebraic extensions of $H$. We briefly review now on these two topics.

A subgroup $H\leqslant F$ is called a \emph{retract of $F$} (just \emph{retract} if there is no risk of confusion) if
there exists a homomorphism $\rho \colon F \to H$ which fixes the elements of $H$ (i.e., such that $\rho^2 =\rho$);
such a morphism is called a \emph{retraction}. The obvious examples of retracts are the free factors of $F$, but there
are retracts which are not free factors. Recognizing retracts is algorithmically possible, as showed
in~\cite[Proposition~4.6]{MVW} following an argument indicated by E.~Turner, though quite complicated in practice,
because it makes use of Makanin's algorithm to solve systems of equations in free groups.

\begin{prop}[4.6 in~\cite{MVW}]\label{retracts}
Let $H\leqslant F$ be a finitely generated subgroup of $F$, given by a finite set of generators. It is algorithmically
decidable whether $H$ is a retract of $F$ and, in the affirmative case, find a retraction $\rho \colon F\to H$.
\end{prop}

For a given endomorphism $\varphi \colon F \to F$, define the \emph{stable image} of $\varphi$ as
$F\varphi^{\infty}=\cap_{m=1}^{\infty} F\varphi^m$. With a simple argument, W.~Imrich and E.~Turner showed in~\cite{IT}
that: 1) $F\varphi^{\infty}$ is a $\varphi$-invariant subgroup of $F$; 2) the restriction of $\varphi$ to its stable
image is always an automorphism; and 3) $F\varphi^{\infty}$ is a retract of $F$. This will be used later in order to
reduce a certain computation with endomorphisms, to a similar one with automorphisms.

\medskip

Let $H\leqslant K\leqslant F$. We say that the extension $H\leqslant K$ is \emph{algebraic}, denoted $H\leqalg K$, if
$H$ is not contained in any proper free factor of $K$. The antagonistic situation consists of $H$ being a \emph{free
factor} of $K$, denoted $H\leqff K$. It is not difficult to see (see~\cite{MVW}) that every extension $H\leqslant K$ of
finitely generated (free) subgroups of $F$ can be decomposed, in a unique way, as an algebraic extension followed by a
free extension, namely $H\leqalg L\leqff K$ (just take $L$ to be the smallest free factor of $K$ containing $H$, or the
biggest algebraic extension of $H$ contained in $K$). The uniqueness refers to the fact that $L$ is completely
determined by the original extension $H\leqslant K$; again, mimicking the classical Galois theory, $L$ is called the
\emph{algebraic closure of $H$ in $K$}. We refer the reader to~\cite{MVW} for a detailed development of these ideas,
including an analysis of the similarities and the significant differences with respect to classical field theory.

The important fact in this story is an old result by M. Takahasi, originally proven by combinatorial methods
in~\cite{T} (and reproduced in Section 2.4, Exercise 8, of~\cite{MKS}). However, the modern graphical techniques
developed by Stalling's in the 1980's (see~\cite{St}) lead to a new, clear, concise and very natural proof of
Takahasi's Theorem, which was discovered independently by E.~Ventura in~\cite{V2}, and by I.~Kapovich and A.~Miasnikov
in~\cite{KM}. S.~Margolis, M.~Sapir and P.~Weil, also independently, considered the same construction in~\cite{MSW} for
a slightly different purposes. See~\cite{MVW} for a unification of these three points of view, written in the language
of algebraic extensions. In this setting, Takahasi's Theorem says the following:

\begin{thm}[Takahasi]\label{tak}
Let $H\leqslant F$ be a subgroup of a free group $F$. If $H$ is finitely generated then it has finitely many algebraic
extensions, i.e.
 $$
\mathcal{AE}(H)=\{ K\leqslant F \mid H\leqalg K\}
 $$
is finite. Furthermore, the elements in $\mathcal{AE}(H)$ are finitely generated, and bases of all of them are
computable from any given set of generators for $H$.
\end{thm}

\noindent {\it Sketch of proof} (see~\cite[Proposition~3.7]{MVW} for details). Think $F =\langle A \mid \,\, \rangle$
as the fundamental group of a bouquet of $n$ circles, and then $H$ as the corresponding covering $X(H)$, which can be
though of as a graph with labels from $A$ on the edges (this graph is easily computable from any given set of
generators of $H$). When $H$ is of infinite index in $F$, the graph $X(H)$ is infinite but, if $H$ is finitely
generated, $X(H)$ consists on a finite core $\Gamma(H)$ with attached infinite trees (each isomorphic to a connected
subgraph of the Cayley graph of $F$ with respect to $A$). Now consider an arbitrary extension $H\leqslant K\leqslant
F$. It corresponds to another covering $X(K)$, which is in turn covered by $X(H)$. That is, $X(K)$ is a quotient of
$X(H)$ and so, can be obtained from $X(H)$ by performing several identifications of vertices and edges. These
identifications may destroy $\Gamma(H)$, but some quotient of $\Gamma(H)$ always remains as a subgraph of $X(K)$ (in
fact, of $\Gamma(K)$). If $H$ is finitely generated then $\Gamma(H)$ is finite, and so has finitely many quotients,
which are computable from $\Gamma(H)$ (i.e. from any given set of generators of $H$). This gives a computable finite
list of extensions of $H$, called the \emph{fringe of $H$}, $\mathcal{O}(H)=\{ H_1, \ldots ,H_p \}$, $p\geqslant 1$.
And, by construction, it is clear that, for every $H\leqslant K$, there exists $i=1,\ldots ,p$ such that $H\leqslant
H_i\leqff K$. This implies that $\mathcal{AE}(H)\subseteq \mathcal{O}(H)$ and so, we already have a proof of the
finiteness part of Takahasi's Theorem. Unfortunately, the equality between these two sets is not true in general, as
one can possibly find free factor relations between the $H_i$'s. But, after a cleaning process (checking for every pair
$(i,j)$ whether $H_i \leqff H_j$ and, in this case, deleting $H_j$ from the list) one can algorithmically compute
$\mathcal{AE}(H)=\{ H_1, \ldots ,H_q \}$, $1\leqslant q\leqslant p$. See~\cite{RVW} for a polynomial time algorithm to
check free factorness. \qed

Note that the smallest and the biggest of the $H_i$'s in $\mathcal{O}(H)$ correspond, respectively, to the quotient
identifying nothing, which gives $H$ itself, and to the quotient identifying all vertices down to a single one, which
gives $\langle A'\rangle \leqff F$, where $A'\subseteq A$ is the set of all letters involved in the generators of $H$.
Note also that the first one belongs to $\mathcal{AE}(H)$ (since $H\leqalg H$) and the same happens for either the
second one or a free factor of it. In particular, $\mathcal{AE}(H)$ contains at least $H$, and a free factor of $F$
(which may coincide). This fact will be used later.

Finally, we mention one of the results in~\cite{St}. Given two finitely generated subgroups $H,K\leqslant F$ (by sets
of generators, say), one can algorithmically compute a basis for $H\cap K$ using the technique of pull-backs of graphs.

\begin{prop}[Stallings,~\cite{St}]\label{pb}
Let $H,K\leqslant F$ be two finitely generated subgroups of a free group $F$, given by finite sets of generators. Then,
a basis for $H\cap K$ is algorithmically computable.
\end{prop}

\section{The algorithms}\label{algo}

Let $H\leqslant F$ be a finitely generated subgroup of $F$, given by a set of generators. We shall give two algorithms
to compute a basis for $\acl_{F} (H)$ and $\ecl_{F} (H)$, respectively. The basic fact that we use is a classical
result due to J.~McCool (see Proposition~5.7 in Chapter~I of~\cite{LS}, and the subsequent paragraph):

\begin{thm}[McCool, \cite{LS}]\label{McCool}
Let $H\leqslant F$ be a finitely generated subgroup of a (finitely generated) free group, given by a finite set of
generators. Then the stabilizer, $\auH(F)$, of $H$ is also finitely generated (in fact finitely presented), and a
finite set of generators (and relations) is algorithmically computable from $H$.
\end{thm}

\subsection{The automorphism case}

By Theorem~\ref{McCool}, $\auH(F)$ is finitely generated; furthermore, a list of generators, say $\auH (F)=\langle
\varphi_1,\ldots ,\varphi_m\rangle \leqslant \au (F)$, can be algorithmically found from a set of generators of $H$.
Now it is clear that
 $$
\acl_{F} (H)=\fix (\auH (F)) =\cap_{\varphi \in \auH(F)} \fix(\varphi) =\fix(\varphi_1)\cap \cdots \cap
\fix(\varphi_m).
 $$
By Maslakova's Theorem~\ref{mas}, we can then compute generators for each of the $\fix(\varphi_i)$'s and, using
Proposition~\ref{pb}, find a basis for their intersection i.e., $\acl_{F} (H)$. Finally, Proposition~\ref{2n} ensures
us that a certain subset of at most $2r(F)$ of those $\varphi_i$'s also makes the job; it only remains to recurrently
compute intersections until finding such a set. Thus, we have proven

\begin{thm}\label{a}
Let $H\leqslant F$ be a finitely generated subgroup of a free group, given by a finite set of generators. Then, a basis
for the auto-closure $\acl_{F} (H)$ of $H$ is algorithmically computable, together with a set of $m\leqslant 2r(F)$
automorphisms $\varphi_1,\ldots ,\varphi_m \in \au (F)$, such that $\acl_{F} (H)=\fix(\varphi_1) \cap \cdots \cap \fix
(\varphi_m)$.\phantom{a} \qed
\end{thm}

\begin{cor}
Let $H\leqslant F$ be a finitely generated subgroup of a free group, given by a finite set of generators. Then, it is
algorithmically decidable whether $H$ is auto-fixed and, in the affirmative case, find a set of $m\leqslant 2r(F)$
automorphisms $\varphi_1,\ldots ,\varphi_m \in \au (F)$, such that $H=\fix(\varphi_1) \cap \cdots \cap \fix
(\varphi_m)$.
\end{cor}

\demo Apply Theorem~\ref{a} to $H$. If $\acl_{F} (H)$ is strictly bigger than $H$, then $H$ is not auto-fixed (there
are elements outside $H$ which are fixed by \emph{every} automorphism of $F$ fixing $H$). Otherwise, $\acl_{F} (H)=H$
and the algorithm in Theorem~\ref{a} also outputs a list of $m\leqslant 2r(F)$ automorphisms $\varphi_1,\ldots
,\varphi_m \in \au (F)$, such that $\fix(\varphi_1) \cap \cdots \cap \fix (\varphi_m) =\acl_{F} (H)=H$. \qed

We don't play much attention to the complexity of this algorithm because it seems far from fast. McCool's algorithm is
a brute force search which is not conceptually complicated, but has strongly exponential complexity. And Maslakova's
algorithm is conceptually much more sophisticated, and its complexity also seems to be quite high. Finally, the
algorithm to compute intersections is both easy and fast.

\subsection{The endomorphism case}

There is no hope that a similar strategy could work in general for endomorphisms instead of automorphisms. On one hand,
Maslakova's Theorem makes strong use of train tracks, a machinery that only works for monomorphisms and definitely does
not work in presence of non-trivial kernel; in fact, at the time of writing, no algorithm is known to compute the fixed
subgroup of an arbitrary endomorphism of $F$. On the other hand, as the following example shows, $\enH (F)$ is not
always finitely generated as submonoid of $\en (F)$ and so, there is no hope of having a variation of McCool's result
for endomorphisms.

\begin{example}[Ciobanu-Dicks,~\cite{CD}]\label{ex} {\rm
We reproduce here example~1.4 of~\cite{CD} to show that $\enH (F)$ is not always finitely generated as a submonoid of
$\en (F)$, even with $H$ being finitely generated as a subgroup of $F$.

Let $F=\langle a,b,c\rangle$ be the free group of rank 3, let $d=ba[c^2,b]a^{-1}$ (where $[x,y]=xyx^{-1}y^{-1}$), and
consider the subgroup $H=\langle a,d\rangle \leqslant F$. Consider the endomorphism $\psi \colon F\to F$ given by
$a\mapsto a$, $b\mapsto d$, $c\mapsto 1$, and the automorphism $\phi \colon F\to F$ given by $a\mapsto a$, $b\mapsto
b$, $c\mapsto cb$. Straightforward computations show that $d\psi =d$ hence, $\psi \in \enH (F)$. Moreover, $\phi^n
\psi$ acts as $a\mapsto a$, $b\mapsto d$, $c\mapsto d^n$ and so, we also have $\phi^n \psi \in \enH (F)$ for every
$n\in \mathbb{Z}$. Now, Corollary~3.4 from~\cite{CD} shows that this is precisely the whole stabilizer of $H$,
 $$
\enH (F) =\{ 1,\, \phi^n\psi \mid n\in \mathbb{Z} \} =\{ 1\} \cup \langle \phi\rangle \psi.
 $$
Again, an easy calculation shows that $(\phi^n\psi )\cdot (\phi^m\psi )=\phi^n\psi$, for every $n,m\in \mathbb{Z}$.
Thus, the monoid $\enH (F)$ is \emph{not} finitely generated. $\Box$}
\end{example}

Being convinced that the above algorithm for the automorphism case cannot be adapted to the endomorphism case, a
different strategy is needed. We shall use algebraic extensions, Takahasi's Theorem and retractions to reduce the
computation of the endo-closure $\ecl_{F} (H)$ to finitely many computations of auto-closures.

\begin{thm}\label{e}
Let $H\leqslant F$ be a finitely generated subgroup of a free group, given by a finite set of generators. Then, a basis
for the endo-closure $\ecl_{F} (H)$ of $H$ is algorithmically computable, together with a set of $m\leqslant 2r(F)$
endomorphisms $\varphi_1,\ldots ,\varphi_m \in \en (F)$, such that $\ecl_{F} (H)=\fix(\varphi_1) \cap \cdots \cap \fix
(\varphi_m)$.
\end{thm}

\demo Consider the set of algebraic extensions of $H$, $\mathcal{AE}(H)=\{ H_1, H_2,\ldots ,H_q \}$, and the subset of
those which are retracts of $F$, say $\mathcal{AE}_{ret}(H)=\{ H_1,\ldots ,H_r \}$, $1\leqslant r\leqslant q$ (note
that $\mathcal{AE}_{ret}(H)$ is not empty because, as we noted above, $\mathcal{AE}(H)$ contains at least a free factor
(and so a retract) of $F$). By Theorem~\ref{tak}, we can algorithmically compute $q\geqslant 1$, and a basis for each
$H_1,\ldots ,H_q$. Now, using Theorem~\ref{retracts}, we can algorithmically decide which of these $H_i$'s are retracts
of $F$, and so compute $r\geqslant 1$, $\mathcal{AE}_{ret}(H)=\{ H_1,\ldots ,H_r \}$, and retractions $\rho_i \colon F
\to H_i$, for $i=1,\ldots ,r$. Then, write the generators of $H$ in terms of the computed bases of each one of these
$H_i$'s, and apply $r$ times Theorem~\ref{a} to compute, for every $i=1,\ldots ,r$, a basis for $\acl_{H_i}(H)$
together with a collection of (at most $2r(H_i)$) automorphisms $\alpha_{i,j} \in \au (H_i)$ such that $\cap_j \fix
(\alpha_{i,j}) =\acl_{H_i} (H)$, and bases for all these fixed subgroups $\fix (\alpha_{i,j})$. Finally, use
Proposition~\ref{pb} to find a basis for $\cap_{i=1}^r \acl_{H_i}(H)$.

Now, we claim that $\cap_{i=1}^r \acl_{H_i}(H)=\ecl_{F}(H)$. In fact, we shall prove this equality under the form
\begin{equation}\label{1}
\bigcap_{i=1}^r \bigcap_{\scriptsize \begin{array}{c}\alpha \in \au(H_i) \\ H\leqslant \fix (\alpha)
\end{array}} \fix (\alpha) =\bigcap_{\scriptsize \begin{array}{c}\beta \in \en(F) \\ H\leqslant \fix (\beta)
\end{array}} \fix (\beta),
\end{equation}
by showing that every intersecting subgroup in one side in also present in the opposite side.

Let $\beta \in \en(F)$ be such that $H\leqslant \fix (\beta)$. Consider the stable image of $\beta$, which is a retract
of $F$, and contains $\fix (\beta)$ and so $H$; then, look at the algebraic closure of $H$ in it,
 $$
H\leqalg H_i \leqff F\beta^{\infty}\leqslant_{\rm ret} F.
 $$
Since free factors of retracts of $F$ are retracts of $F$, this $H_i$ is an element of $\mathcal{AE}_{ret}(H)$.
Furthermore, the endomorphism $\beta$ restricts to an automorphism of $F\beta^{\infty}$ which, in turn, restricts to an
automorphism $\alpha=\beta_{|H_i}$ of $H_i$ (because images of free factors of $F\beta^{\infty}$ under $\beta$, are
again free factors of $F\beta^{\infty}$). And, clearly, $H\leqslant \fix (\alpha) \leqslant \fix (\beta)$. This shows
inclusion ``$\leqslant$" in equation~(\ref{1}).

Now, let $H_i\in \mathcal{AE}_{ret}(H)$, and let $\alpha \in \au(H_i)$ with $H\leqslant \fix (\alpha)$. Consider $\beta
=\rho_i \alpha \iota_i \in \en (F)$, where $\rho_i \colon F\to H_i$ is a retraction, and $\iota_i \colon H_i \to F$ is
the inclusion map. It is clear that $H\leqslant \fix (\alpha)=\fix (\beta)$ and so, $\fix(\alpha )$ is also one of the
subgroups appearing the intersection on the right hand side of~(\ref{1}). This shows inclusion ``$\geqslant$" in
equation~(\ref{1}) and completes the proof of the claim.

Thus, the algorithm described in the first paragraph of this proof, certainly computes a basis for $\ecl_{F}(H)$
(together with some side information, namely the retractions $\rho_i$, the collection of automorphisms $\alpha_{i,j}
\in \au (H_i)$, and bases for their fixed subgroups $\fix (\alpha_{i,j})$). To conclude, it only remains to explicitly
construct a list of at most $2r(F)$ endomorphisms of $F$ whose fixed set is exactly $\ecl_{F}(H)$. This is easy from
the previous paragraph: the collection of endomorphisms $\beta_{i,j} =\rho_i \alpha_{i,j} \iota_i \in \en (F)$ satisfy
$\fix (\beta_{i,j}) =\fix (\alpha_{i,j})$ and so,
 $$
\bigcap_{i,j} \fix (\beta_{i,j}) = \bigcap_i \big( \bigcap_j \fix (\alpha_{i,j}) \big) = \bigcap_i \acl_{H_i} (H)
=\ecl_F (H).
 $$
It could happen that the computed set, $\{ \beta_{i,j} \mid i,j \}$, of endomorphisms of $F$ exceeded in number the
maximum wanted quantity of $2r(F)$. In this case, Proposition~\ref{2n} ensures us that a certain subset of cardinal at
most $2r(F)$ makes the job, too. Since, as a side product of our computation, we also have a basis of each $\fix
(\beta_{i,j}) =\fix (\alpha_{i,j})$, it only remains to recurrently compute intersections until finding the desired set
(knowing it exists). This concludes the proof. \qed

\begin{cor}\label{endofixedness}
Let $H\leqslant F$ be a finitely generated subgroup of a free group, given by a finite set of generators. Then, it is
algorithmically decidable whether $H$ is endo-fixed and, in the affirmative case, find a set of $m\leqslant 2r(F)$
endomorphisms $\varphi_1,\ldots ,\varphi_m \in \en (F)$ such that $H=\fix(\varphi_1) \cap \cdots \cap \fix
(\varphi_m)$. \qed
\end{cor}

\section{Open problems}\label{p}

In this last section we collect a list of interesting questions and open problems in this subject.

\medskip

\noindent \textbf{Problem 1.} Find an algorithm to compute $\fix (\varphi)$ for a given $\varphi \in \en (F)$.

\medskip

\noindent \textbf{Problem 2.} Find an algorithm to determine whether a given finitely generated subgroup $H\leqslant F$
is 1-auto-fixed, or 1-endo-fixed.

\medskip

\noindent \textbf{Problem 3.} Do the families of auto-fixed and 1-auto-fixed subgroups of $F$ coincide ? And those of
1-endo-fixed and endo-fixed subgroups ?

\medskip

\noindent \textbf{Problem 4.} Find effective, say polynomial time, algorithms to compute the fix subgroup of a given
endomorphism, and to determine whether a given finitely generated subgroup $H$ of $F$ is 1-auto-fixed, or 1-endo-fixed,
or auto-fixed, or endo-fixed.

\medskip

\noindent \textbf{Problem 5.} Are endo-fixed subgroups of $F$ $F$-inert ?

\medskip

\noindent \textbf{Problem 6.} Find an algorithm to decide whether a given subgroup $H\leqslant F$ is $F$-inert.

\end{document}